\documentclass{amsart}

\begin{document}

\newtheorem{thm}{Theorem}[section]
\newtheorem{lem}[thm]{Lemma}
\newtheorem{cor}[thm]{Corollary}

\theoremstyle{definition}
\newtheorem{defn}{Definition}[section]

\theoremstyle{remark}
\newtheorem{rmk}{Remark}[section]

\def\square{\hfill${\vcenter{\vbox{\hrule height.4pt \hbox{\vrule
width.4pt height7pt \kern7pt \vrule width.4pt} \hrule height.4pt}}}$}
\def\T{\mathcal T}

\newenvironment{pf}{{\it Proof:}\quad}{\square \vskip 12pt}

\title{Foliations Transverse to Triangulations of $3$-Manifolds}
\author{Danny Calegari}
\address{Department of Mathematics \\ UC Berkeley \\ Berkeley, CA 94704}
\email{dannyc@math.berkeley.edu}

\maketitle

\begin{abstract}
We investigate the combinatorial analogues, in the context of normal
surfaces, of taut and transversely measured (codimension $1$) foliations
of $3$-manifolds. We establish that the existence of certain combinatorial
structures, a priori weaker than the existence of the corresponding
foliation, is sufficient to guarantee that the manifold in question 
satisfies certain properties, e.g. irreducibility. The finiteness of our
combinatorial structures allows us to make our results quantitative in nature
and has (coarse) geometrical consequences for the manifold. Furthermore, our
techniques give a straightforward combinatorial proof of Novikov's theorem.
\end{abstract}

\section{Introduction}

In this paper we study some of the relationships between normal surfaces and
foliations. Our general approach is to discuss what combinatorial structure
on a triangulation of a $3$-manifold is sufficient to guarantee the existence
of a foliation in ``normal form'' with respect to the triangulation. i.e.
we want every leaf of the foliation to be a normal surface away from the
vertices of $M$. Or, if we cannot guarantee such a foliation, we consider
what combinatorial structure is sufficient to guarantee properties of $M$
that would be a consequence of the existence of a foliation.

The combinatorial structure we describe consists of a choice of orientation
for each edge of the $1$-skeleton $M^1$. That there exists a transverse
foliation locally amounts to a local condition on the star of each vertex.
The problem of the global existence of a foliation seems to be a very
hard problem, and we are not able to treat it effectively except in some
special cases. In particular, the related problem of when a branched surface
carries (abstractly) a lamination was shown to be algorithmically unsolvable,
by Lee Mosher (see {\bf [Ga]}), although it is not known whether the problem is
still unsolvable if the branched surface is given together with an
embedding in a $3$-manifold.

\vskip 12pt

In the first section, we are able to give precise conditions for the existence
of a {\em transversely measured} normal foliation. Our condition says 
that directed loops in the $1$-skeleton (with respect to our orientation)
must lie in an open half-space of homology. As a corollary, we are able to
give an elementary combinatorial proof that a $3$-manifold admitting a taut
foliation in which the transverse loops lie in an open half-space of homology,
is a surface bundle over $S^1$.

\vskip 12pt

In the second section we treat the converse question, of when a foliation of
$M$ by closed surfaces can be put in normal form with respect to a 
triangulation. We show that if one can make the $1$-skeleton transverse
to the foliation and have at least one incoming and one outgoing edge from
each vertex, then the foliation can be isotoped rel. its intersection with
the $1$-skeleton to be in normal form.

\vskip 12pt

In the next two sections we discuss when a triangulation admits a taut
foliation in normal form. If a choice of orientations on the $1$-skeleton
as above admits a transverse foliation locally, if every oriented loop is
homotopically essential, and if in addition the $1$-skeleton is recurrent,
thought of as a directed graph, then we are able to show that the universal
cover of $M$ can be given a foliation in normal form such that each leaf
is incompressible. As a corollary, we see that a manifold admitting such a
combinatorial structure is either irreducible or $S^2\times S^1$. This
is a generalization of Novikov's theorem, and it remains to be seen whether
the combinatorial structure so described is more general than that of the
existence of a taut foliation. Our methods give new information even about
$3$-manifolds admitting taut foliations, showing that any deformation of
the foliation in the universal cover which is ``roughly equivariant'' (i.e.
preserves the normal disk types in each tetrahedron) also has incompressible
leaves. Moreover, our technique gives an elementary (combinatorial) proof
that circles transverse to leaves of a taut foliation
are homotopically essential. 

In fact, we can make our results quantitative, and show that the leaves in
the universal cover satisfy the same kind of isoperimetric inequalities 
that minimal surfaces in $M$ satisfy. Our argument here uses the finiteness
of the combinatorial structure, which makes uniform geometric estimates 
almost effortless.

\vskip 12pt

In the penultimate section, we weaken our hypotheses considerably and show
that we can still get a strong result. In particular, we weaken the
condition that oriented loops be homotopically essential to
the condition that oriented loops of length $\le k_1$ bound no simplicial
disks of simplicial area $\le k_2$, for constants $k_1,k_2$
depending on the triangulation and the choice of orientation. With this
hypothesis, we can nevertheless show that $\pi_2(M)=1$ or $M=S^2\times S^1$.

This result is interesting insofar that it shows that even the ``local"  
existence of a foliation transverse to the $1$-skeleton is enough to
get topological information about the manifold.

In particular, we
expect that this condition is much more general than the existence of a
foliation in normal form on $M$ compatible with the orientations.
 This brings to mind
a question of D. Gabai in {\bf [Ga]}, namely: ``Do there exist useful
branched surfaces which do not carry anything''? Though we do not demonstrate
the non-existence of foliations compatible with our combinatorial structures,
nevertheless, we do not need to produce a globally compatible foliation to find
the combinatorial structure ``useful''.

\vskip 12pt

We adhere to the convention, it what follows, even if we do not mention it
explicitly, that every foliation is oriented, co-oriented, and smooth. For
simplicity, and to avoid headaches, we have not investigated the extent to
which these conditions can be dropped.

\vskip 12pt

I would like to thank Andrew Casson for his patience, his comments, and his
suggestions regarding the following material. I am also grateful to the
referee for some excellent observations and comments.

\section{Triangulations and maps to $S^1$}

Let $M$ be a compact $3$-manifold, and let $\mathcal{T}$ 
be a triangulation of $M$. 
A {\it positive orientation} on $\mathcal{T}$ 
is a choice of orientation on each edge of
$\mathcal{T}^1$, and a choice of $\alpha \in H^1(M;\mathbb{R})$ such that
\begin{itemize}
\item{For each loop $\gamma \subset \mathcal{T}^1$ 
with $\gamma$ positively oriented,
$\alpha([\gamma])>0$, where $[\gamma]$ denotes the image of $\gamma$ in 
$H_1(M;\mathbb{R})$}
\item{For every vertex $v \in \mathcal{T}^1$ 
let $o(v)$ be the set of vertices $w$
such that there is an oriented edge from $v$ to $w$, and let $i(v)$ be the
set of vertices $w$ such that there is an oriented edge from $w$ to $v$. 
Then there
is a connected subgraph of ${\rm link}(v)$ whose vertices are exactly $o(v)$, 
and a connected subgraph of ${\rm link}(v)$ whose vertices are exactly $i(v)$.}
\end{itemize}

We call a choice of orientation for each edge of a cell complex $C$ a
{\it direction} on $C$. Maps between directed $CW$ complexes 
are orientation-preserving
if they preserve the orientation on each edge.

\begin{lem}
Let $C$ be a complex together with a direction.
There is an orientation-preserving immersion $f:C^1 \to S^1$ 
which is affine on each edge (with respect to the standard affine structure
on $S^1$), and extends to
all of $C$ iff all oriented cycles in $C^1$, considered as elements of
$H_1(C)$, are contained in an open half-space.
\end{lem}
\begin{pf}
An orientation-preserving immersion $f:C^1 \to S^1$ is determined by
a collection of real values $x_i>0$, one for each edge $e_i$, such that we have
$$\sum_{e_i \in \gamma} x_i \in \mathbb{Z}$$
for each oriented cycle $\gamma$. This map extends to $C^2$ iff for each
disk $D \in C^2$, the boundary $\partial D$ 
is mapped to $S^1$ null-homotopically.
That is, writing the boundary $\partial D$ as a union of positively
and negatively oriented edges,
$$E_j = \sum_{e_i \in \partial D_j} x_i - 
\sum_{-e_i \in \partial D_j} x_i = 0$$
where our notation is meant to indicate whether the orientation of $e_i$ agrees
or disagrees with the orientation of $\partial D$ induced by some arbitrary
orientation on $D$.

It is a fact that the collection of equalities $E_j = 0$ has a solution 
with all
$x_i$ positive iff there is no linear combination of $\sum_j v_j E_j$ such that
$$\sum_j v_j E_j = \sum_i c_i x_i$$ with all $c_i$ non-negative, 
and at least one
$c_i$ positive. For, let $H$ be the subspace of $\mathbb{R}^n$ spanned by
the vectors $E_i$. Then the orthogonal subspace $H^\perp$ is precisely the
set of solutions to the equations $E_i = 0$. If $K$ denotes the convex cone
where all $x_i>0$, then vectors $v\in K \cup -K$ are characterised by the
property that $v^\perp$ does not intersect $K$. From certain perspectives,
this is just the finite dimensional version of the Hahn-Banach theorem. 

Notice that, 
thought of as an element of the set of $1$-chains on $M$, each $E_j$
is a boundary. Therefore $\sum_j v_j E_j = \partial a$ for some $2$-chain $a$.
But any sum $\sum_i c_i x_i$ with all $c_i$ non-negative that represents a
$1$-cycle must consist of a non-negative sum of positively oriented cycles 
$\gamma \in C^1$. Since such cycles all lie in the same open half-space of
$H_1$, a positive combination of them cannot be a boundary.

Having found a positive solution, 
we may find one with rational coefficients, and
therefore by scaling, with integral coefficients.

Finally, since $S^1$ is aspherical, any map on $C^2$ extends to $C$.
\end{pf}

With this lemma we may establish

\begin{thm}
If $M$ admits a triangulation with a positive orientation, then $M$ is 
a surface bundle over $S^1$ and there is a projection map $M \to S^1$ 
which is affine on each simplex.
\end{thm}
\begin{pf}
Let $\mathcal{M}$ denote a geometric model for $M$ constructed by letting each
tetrahedron of the triangulation be a regular Euclidean tetrahedron with side
lengths equal to $1$.

By our lemma, the first condition implies that there
is an orientation-preserving immersion 
$f:\Gamma \to S^1$ which extends to all of
$\mathcal{T}$.
By a genericity assumption, we insist that the images of all the vertices
are sent to distinct points of $S^1$.
The restriction of $f$ to the $1$-skeleton of each tetrahedron $\Delta$ 
therefore
lifts to a map $\tilde f:\Delta^1 \to \mathbb{R}$. We extend $\tilde f$ to
the entire tetrahedron $\Delta$ by requiring it to be affine.
It is clear that this can be done
compatibly to give an piecewise-affine map $\mathcal{M} \to S^1$. We claim that
this map is the projection map from the total space of a surface bundle to the
base space. To see this, it suffices to show that the foliation 
of $\mathcal{M}$ by the preimages of points in $S^1$ is non-singular.

Notice that there is only one possible orientation for each 
$\Delta^1$, up to isomorphism, 
and it is clear that this orientation induces a non-singular foliation
on $\Delta$. This foliation pieces together compatibly along faces and edges.
It remains to check that it is non-singular at vertices.

Let $v$ be a vertex of the triangulation, and ${\rm star}(v)$ denote the
star of $v$. Then $f$ lifts to
$\tilde f: {\rm star}(v) \to \mathbb{R}$, since topologically, 
${\rm star}(v)$ is a $B^3$ which is
simply-connected. ${\rm link}(v)$ is a triangulated $S^2$. 
Then $f|_{{\rm link}(v)}$
is non-degenerate away from the vertices, 
and its level sets away from these consist
of a disjoint union of circles.

If $\lambda$ is the
leaf containing $v$, then $\lambda \cap {\rm star}(v)$ is the cone on the 
set of
$p \in {\rm link}(v)$ where $\tilde f(p) = \tilde f(v)$. We must show that this
is a circle. By genericity, this set is a disjoint union of circles. 
Let $S$ be a
circle in ${\rm link}(v)$ separating a maximal graph whose vertices are $o(v)$ 
from a maximal graph whose vertices are $i(v)$. 
Call these two graphs $\Gamma_o$ and
$\Gamma_i$. Such a circle exists, since the
two graphs are disjoint, 
and any two connected closed subsets of $S^2$ are separated
by an embedded circle. 

The union of simplices in ${\rm link}(v)$ intersecting $S$ non-trivially is an
open annulus whose boundary is contained in $\Gamma_i \cup \Gamma_o$,
and it is clear that we can write the annulus as $S^1 \times I$
where $\tilde f$ is monotonically increasing on each $p \times I$. Since the
value of $\tilde f$ on one boundary component of 
this annulus is strictly greater
than $\tilde f(v)$, and strictly less than it on the other component, 
$\tilde f^{-1}(\tilde f(v))\cap S^1 \times I$ is a single circle. Now
for any $p$ in the complement of this annulus, $\tilde f(p)$ is a convex
combination of values strictly larger than, or a convex combination of values
strictly smaller than $\tilde f(v)$ --- namely the values of $\tilde f$ on the
vertices on the appropriate side. Hence
$\tilde f^{-1}(\tilde f(v))\cap {\rm link}(v)$ 
is a single circle, and the foliation
is nonsingular at $v$.
\end{pf}

\begin{cor}
Let $M$ be a $3$-manifold with a taut oriented, co-oriented, smooth foliation 
$\mathcal{F}$ and
let $\alpha \in H^1(M;\mathbb{R})$ be 
such that for every transverse, 
positively oriented cycle $\gamma$ the inequality
$\alpha([\gamma])>0$ is true. Then $M$ is a surface bundle over $S^1$.
\end{cor}

\begin{pf}
To see that $M$ fibers over $S^1$, just take a very fine triangulation of $M$
with all edges transverse to $\mathcal{F}$. Orient the edges according to the
co-orientation on $\mathcal{F}$, and make the edges ``straight" enough that
the induced orientation on the $1$-skeleton of the triangulation is a positive
triangulation.

Essentially, if we restrict to a small $I^3$ foliated by
$z=const.$, we want the triangulation in that small $I^3$ to be by
approximately affine tetrahedra with respect to the affine structure on $I^3$.
The condition is automatically satisfied for affine triangulations, with
edges oriented by a co-orientation on an affine foliation.

For, let $p,q$ be two vertices in $o(v)$. 
Let $\pi$ be the plane spanned by $p,q,v$.
This plane intersects ${\rm link}(v)$ in a circle, since ${\rm star}(v)$ 
is star-shaped
with center $v$. Let $\gamma$ be the arc of this circle joining $p$ to $q$
such that the value of $z$ on $\gamma$ is greater than the value of $z(v)$.
Then either $\gamma$ is an edge of the $1$-skeleton, in which case
$\Gamma_o$ connects $p$ to $q$, or it intersects the $1$-skeleton at some
first point $m$ in the interior of an edge $e$ If the latter case, $z(m)>z(v)$,
so there is a vertex $r$ on $e$ with $z(r)>z(m)>z(v)$. 
We ``slide" $\gamma$ along
$m$ to $r$, and continue inductively to produce a path in $\Gamma_o$ connecting
$p$ to $q$, so $\Gamma_o$ is connected. Similarly, $\Gamma_i$ is connected.

Since Riemannian manifolds are locally almost affine, 
this can be done compatibly over
the entire triangulation. More precisely, since $\mathcal{F}$ is smooth,
we can cover $M$ with co-ordinate patches such that in each co-ordinate
patch, $\mathcal{F}$ is a foliation of $\mathbb{R}^3$ by level sets of
the form $z = {\rm const.}$ and the co-ordinate transformations are ``almost"
linear. Here ``almost" means sufficiently close that the triangulation
can be straightened to an affine triangulation in each patch without
disturbing the local combinatorial structure. Essentially, we just need to
pick a triangulation by sufficiently ``squat'' simplices. A rigorous proof
of this fact can be found in {\bf [Be]}, where it is attributed originally
to Thurston.

Since a smooth foliation
locally resembles such an affine foliation to first order, this can be
done compatibly over the entire manifold.

The triangulation of $M$ is therefore positively oriented, and $M$ is a surface
bundle over $S^1$, as required.
\end{pf}

\begin{rmk}
In fact, $\mathcal{F}$ as above carries a transverse measure $\mu$ such that
$\mu(\gamma) = \alpha([\gamma])$ for any transverse, positively oriented cycle
$\gamma$. Suppose, for example, that some leaf $\lambda \in \mathcal{F}$ is
dense in $M$. Then for a given transversal $t$, 
there are points in $t \cap \lambda$
arbitrarily near the endpoints of $t$. 
These can be joined up by a path in $\lambda$
to give a cycle. Evaluating $\alpha$ on this cycle, 
and taking the supremum over
all pairs of points which tend toward the ends of $t$, we get $\mu(t)$. This is
positive, since it is greater than the value of $\alpha$ on some positive
transverse cycle, which is $>0$. 
A choice of a different path in $\lambda$ might
give a different cohomology class, but these would differ by a cohomology class
carried by a loop in $\lambda$, and $\alpha$ 
evaluated on this class is necessarily
$0$. If we homotope $t$, keeping its endpoints on the same leaf, we can join it
up by a path in $\lambda$ homotopic to 
the original path, and therefore giving the
same value when $\alpha$ is evaluated on it. 
Finally, if we write $t$ as a union of
two intervals, then a cycle joining up $t$ is 
homologous to a sum of two cycles,
each joining up one of the sub-intervals of $t$.

If we pick some fine triangulation, 
and an associated map to $S^1$, the pullback of
the angular measure on $S^1$ to $M$ is ``approximately" a 
transverse measure for
$\mathcal{F}$. It is positive on any monotone path in the $1$-skeleton of the
triangulation. However, monotone paths which are ``nearly horizontal" are not
generally approximated by a monotone path in the $1$-skeleton. 
By including more and more monotone paths in the $1$-skeleton, we can
approximate $\mu$ more and more closely. 
In fact, the more elements of $H_1$ that
are carried by monotone cycles in the triangulation, 
the less flexibility we have
in choosing the homotopy class of our map to $S^1$, and the better the pullback
measure approximates $\mu$.
\end{rmk}

\begin{rmk}
Notice that there is nothing inherently $3$-dimensional about our theorem.
If $M$ is an arbitrary manifold, and $\mathcal{F}$ a smooth codimension one
co-oriented foliation such that every transverse cycle lies in a half-space
of $H_1$, then if we choose a very fine triangulation transverse to the
foliation which can be straightened to an affine triangulation in every 
co-ordinate chart without disturbing the combinatorial structure, and we
orient the $1$-skeleton according to $\mathcal{F}$, then the piecewise
affine map to $S^1$ guaranteed by the lemma induces a non-singular foliation
everywhere, and exhibits the manifold as a bundle over $S^1$. We prove our
theorem in the $3$-dimensional case only because that is our interest for
applications.
\end{rmk}

\begin{rmk}
In {\bf [Su]}, D. Sullivan proves the following conjecture of R. Edwards:
if $M$ is a foliated manifold with all leaves compact, such that the homology
classes represented by the leaves lie in an open half-space of homology, for
the appropriate dimension, then $M$ is transversely measured. Notice that
if the foliation is of codimension $1$, then this condition implies the
condition of the corollary above. Our proof owes something to Sullivan's
approach - in particular, the key result from linear algebra that we use is
a finite dimensional version of the Hahn-Banach theorem, which in its
general form is essential in setting up the machinery for Sullivan's theorem.
\end{rmk}

\section{Normal Form for Surface Bundles}

The results of the last section suggest the question of when a triangulation
of a surface bundle admits a positive orientation.

\begin{thm}
Let $M$ be a surface bundle over $S^1$ with
projection map $r:M \to S^1$. Let $\mathcal{T}$ be any triangulation of
$M$ such that the star of each vertex is an embedded $B^3$ in $M$, and
such that the following two conditions are satisfied:

\begin{itemize}
\item{each edge is monotone with respect to $r$ and oriented according to
the orientation in the circle direction}
\item{there is an outgoing edge and an incoming edge for every vertex}
\end{itemize}

Assume that the genus of the surface is at least $1$.

Then the induced orientation of the edges is a positive orientation.
\end{thm}
\begin{pf}

We lift this triangulation to a triangulation of the universal cover, which
is $\mathbb{R}^3$ foliated by planes. Notice that the lift of the triangulation
also satisfies these two properties. By abuse of notation, we denote by
$r$ the map $r:\tilde M \to \mathbb{R}$ whose preimages are the lifts of
surfaces.

By the results of the first section, there is an (equivariant) map
$f:\tilde M \to \mathbb{R}$ affine on each simplex, and agreeing with
$r$ on the $1$-skeleton.

We can further assume a non-degeneracy condition, namely that for each
pair of vertices $v,w$, $f(v)\ne f(w)$.

Then $$C_1 \cup C_2 \cup \dots \cup C_i = f^{-1}(f(v)) \cap {\rm link}(v)$$
is a non-empty disjoint collection of circles. Suppose there are at least
two circles, $C_1,C_2$.

The union of the simplices intersecting $C_i$ is an open annulus, and
we label its boundary components $u_i$ and $l_i$, where $f(l_i)<f(C_i)<f(u_i)$
for each $i$. It is possible that $u_1\cap u_2$ or $l_1 \cap l_2$ are
non-empty, but they cannot both be non-empty, for otherwise $C_1$ would be
a non-separating circle in ${\rm link}(v)$, an absurdity. WLOG, say that
$u_1 \cap u_2$ is empty. Let $U_i$ be the disk consisting of the region
of ${\rm link}(v)$ bounded by $l_i$ and containing $C_i$. We assume $U_1$ and
$U_2$ are disjoint, for otherwise there is an annulus between $C_1$ and
$C_2$ and a circle $C_i$ a meridian of this annulus, and we can replace
one of $C_1,C_2$ by $C_i$ if necessary so that this is satisfied. 

Let $x_i$ be the vertex attaining the
maximum value of $f$ on $U_i$. Since $l_i$ is an unknotted circle, and
since $f(l_i)=r(l_i)<r(v)=f(v)$ on this circle, there is an embedded disk
$D_i$ with $r(D_i)<r(v)$ whose boundary is $l_i$. Since $D_i$ and $U_i$
are embedded with a common boundary, we can arrange that their intersection
is a collection of circles. We perform disk exchanges on these circles to
produce an embedded surface $S_i$ made from pieces of $U_i$ and $D_i$, and
containing ${\rm nbhd}(x_i)\cap U_i$. (In fact, we do not need to do these
exchanges - we can simply take as $S_i$ the boundary of some complementary
region containing $x_i$.)

By assumption, there is an infinite increasing ray contained 
in the $1$-skeleton
emanating from $x_i$. Call this ray $\alpha_i$. Since the value of
$f$ on $\alpha_i$ is greater than $f(x_i)$, this ray cannot intersect
$S_i$ except at $x_i$. For, if it does so, it intersects $S_i$ in the 
$1$-skeleton. But for every point on $S_i$, either $r$ or $f$ is less than
$f(x_i)$, and therefore for the intersection of $S_i$ with the $1$-skeleton,
$f$ is less than $f(x_i)$. Hence $\alpha$ is entirely contained in some
complementary region of $S_i$, and since it is infinite, this region is
unbounded.

Moreover, there is an
increasing edge $e_i$ from $v$ to $x_i$. This edge intersects $S_i$ only
at $x_i$. For, it cannot intersect any piece of $D_i$, since the value of
$r$ is less than $f(v)$ there. Also, it cannot intersect $U_i$ except at
$x_i$, since $e$ intersects ${\rm link}(v)$ only at $x_i$. Therefore the union
$\beta_i = e_i \cup \alpha_i$ is an infinite increasing ray which
intersects $S_i$ exactly once. Hence $v$ is in the bounded complementary
region of $S_i$. 

Notice also that $e_1$ does not intersect $S_2$ at all, nor does $e_2$
intersect $S_1$, since again it can only intersect it in pieces of
$D_i$, and there the value of $r$ is less than $f(v)=r(v)$. Therefore
$x_1$ is contained in the bounded complementary region of $S_2$, and
$x_2$ in the bounded complementary region of $S_1$. WLOG, $f(x_1)>f(x_2)$.
But then $\alpha_1$ cannot intersect $S_2$, since on the intersection of
$S_2$ with the $1$-skeleton, $x_2$ attains the highest value of $f$.
Hence $\alpha_1$ is bounded by $S_2$, which is a contradiction.

Therefore there is exactly one circle $C_1$, and the leaf 
$f^{-1}(f(v))\cap {\rm star}(v)$ is a non-singular disk, and the orientation on
the edges is a positive orientation, as required.

\end{pf}

\begin{rmk}
We may think of this theorem as giving a kind of ``normal form" for surface
bundles with respect to a triangulation. In particular, if the bundle
can be made ``normal" with respect to the $1$-skeleton, this theorem 
guarantees it can be made ``normal" with respect to the entire triangulation. 
\end{rmk}

\begin{rmk}
The issue is to decide for what triangulations $\mathcal{T}$ of a surface
bundle $M$ the $1$-skeleton can be made transverse to the foliation by
surfaces in such a way that each vertex has an outgoing and an incoming
edge, with respect to some co-orientation on the foliation.

Let $r:M \to S^1$ be any map generic with respect to the $1$-skeleton of
$\mathcal{T}$. Then there is a subdivision of the $1$-skeleton to a finite
graph $\Gamma$ such that each edge of $\Gamma$ is transverse to foliation
(i.e. take as additional vertices of $\Gamma$ the critical points of 
$r|_{\mathcal{T}^1}$). There are two issues to be resolved. The first is
whether $r:\Gamma \to S^1$ is homotopic to a map monotone on each edge, and
with an outgoing and an incoming edge from each vertex. The second is
the issue of whether such a homotopy of $r$ can be realized by an isotopy of
$\Gamma$ in $M$.

Such a homotopy can be decomposed into a collection of ``local" moves, which
consist of exchanging the order of neighboring vertices and re-orienting
any edge between the two of them. Let $v_1,v_2$ be the vertices in question,
and suppose $r(v_1)<r(v_2)$. Let $\alpha$ be the segment of $S^1$ between
$r(v_1)$ and $r(v_2)$ containing the image of no other vertex. Then
$r^{-1}(\alpha)$ is homeomorphic to ${\rm surface} \times I$, and
$\Gamma \cap r^{-1}(\alpha)$ consists of a collection of monotone arcs
from one boundary component to the other, together with the set of
outgoing arcs rooted at $v_1$, and the set of incoming arcs rooted at $v_2$,
which are points on opposite boundaries. It is clear that the only obstruction
to performing an isotopy exchanging the order of $v_1$ and $v_2$ is whether
or not the outgoing edges from $v_1$ ``link" the incoming edges to $v_2$.

\end{rmk}

\section{Partial Orderings}

Let $M$ be a closed $3$-manifold, and let $\T$ be a triangulation of
$M$. A {\it direction} on $M$ is a choice of orientation for each edge
in the $1$-skeleton $\T^1$ of the triangulation. A direction is a
{\it local orientation} if it satisfies the conditions

\begin{enumerate}
\item{for each vertex $v$ the maximal subgraphs $o(v)$ and $i(v)$ of
${\rm link}(v)$ whose vertices are, respectively, the outgoing and the incoming
vertices from and to $v$, are nonempty and connected}
\item{the direction restricts to a total ordering on the vertices of each
tetrahedron}
\item{the $1$-skeleton is recurrent as a directed graph. That is, there is
an increasing path from each vertex to each other vertex.}
\end{enumerate}

\vskip 12pt

{\bf Example:} On $S^3$, consider the Hopf vector field. This is a volume
preserving flow, so any cone field which supports this vector field is
recurrent. (for the definition of cone fields, see {\bf [Su]}) If we take
some sufficiently fine triangulation supported by such a cone field, the
local orientation conditions will be satisfied, since locally there is
a product structure given by the flow, which is transverse to our 
triangulation. Again, if the triangulation is sufficiently fine, it can be
made recurrent, since the cone field is recurrent. However, there is no
foliation transverse to this local orientation, for such a foliation would be
taut by recurrence, which is impossible on $S^3$.

\vskip 12pt

Since the vertices of each triangle in $\T^2$ are totally ordered, we
can speak unambiguously of the {\it long} edge of any triangle, and also
of the {\it upper} and {\it lower} short edges. We construct a directed 
graph $\Gamma$ associated to the direction whose vertices are edges of
$\T^1$ and whose directed edges are the ordered pairs $(e_i,e_j)$
where $e_j$ is the long edge, and $e_i$ the upper or lower short edge of
some triangle in $\T^2$. $\Gamma$ is {\it expanding} if it contains
a pair of directed loops, one containing the edge $(e_i,e_k)$, one containing
the edge $(e_j,e_k)$ where $e_i,e_j$ are the upper and lower short, and
$e_k$ the long edge of some triangle.

A choice of orientation on the edges of a triangulation, or more generally
a choice of orientation for the edges of a graph, determines a partial 
ordering on the vertices 
by declaring that $x \le y$ iff there is an
oriented path in the graph or $1$-skeleton from $x$ to $y$. The partial
orderings for the $1$-skeleta of our compact manifolds $M$ will generally
not be very interesting: recurrence implies that for any two elements $x,y$
both $x \le y$ and $y \le x$. However, if we pull back these orientations
to the universal cover of $M$, the induced partial orderings are more 
interesting. We pursue this more vigorously in the next section.

\begin{lem}
If $M$ admits a local orientation in which every oriented loop is homotopically
essential, and if $\Gamma$ is expanding, then $\pi_1(M)$ has exponential 
growth.
\end{lem}
\begin{pf}
By the hypothesis, $\tilde M$ has no oriented loops in its $1$-skeleton.
But then the associated graph $\tilde \Gamma$ has no directed loops. Hence
there is an infinite dyadic tree (the lift of the directed loops guaranteed
by the condition that $\Gamma$ is expanding) which embeds in $\tilde \Gamma$.
Hence $\tilde M$, and therefore $\pi_1(M)$, has exponential growth by the
usual reason that the Cayley graph of $\pi_1(M)$ has the quasi-isometry type
of $M$. (See for instance {\bf [Gr]}).
\end{pf}
\vskip 12pt

The reason to introduce these definitions is given by the following lemma:

\begin{lem}
A co-oriented tautly foliated $3$-manifold admits a triangulation with a 
local orientation in which every oriented loop in the $1$-skeleton is
homotopically essential. Conversely, a foliation in normal form relative
to a local orientation is taut.
\end{lem}
\begin{pf}
As before, choose a triangulation such that $\mathcal{F}$ is in normal
form with respect to the triangulation. Orient the edges of the triangulation
according to the co-orientation on $\mathcal{F}$.
This triangulation can be refined
repeatedly until $M^1$ is recurrent.

The second statement is immediate.
\end{pf}

\begin{thm}
If $M$ admits a local orientation, then the induced orientation on any
connected finite cover of $M$ is a local orientation.
\end{thm}
\begin{pf}
The only non-trivial condition to check is recurrence. Suppose there exists
a monotone path from $p$ to $q$ in the cover. Then this projects to a
monotone path in $M$ which can be completed to a monotone loop, by recurrence.
Then some power of this loop lifts to the cover, so there is a monotone
path from $q$ to $p$. Hence $M^1$ breaks up into recurrent components.
Since it is connected, there is only one such component.
\end{pf}

\begin{thm}
If $M$ admits a local orientation such that every oriented loop is
homotopically essential, and $M$ has a finite cover with fundamental
group $\mathbb{Z}$, then $M= S^2 \times S^1$.
\end{thm}
\begin{pf}
The finite cover, call it $N$, also has the property that every oriented
loop is homotopically essential, and therefore equal to some power of 
the generator
of $H_1(N)$. Suppose there are two oriented loops which represent
$\alpha^n$ and $\alpha^{-m}$ for some positive integers $n,m$. Then there
is another oriented loop representing $\alpha^r$ which connects two points
on these loops. By composing loops, we can find an oriented loop representing
$\alpha^{an - bn + cr}$ for any positive integers $a,b,c$. But this implies
that we can find an oriented loop representing the trivial element, a 
contradiction. Hence all oriented loops lie in an open half-space of
$H_1(N)$, and by our theorem, $N$ is a surface bundle over $S^1$. Since
$\pi_1(N) = \mathbb{Z}$, $N=S^2 \times S^1$. Project a normal $S^2$ down to
$M$. Then the edge weights determined by the (possibly immersed) image 
represent an embedded (possibly disconnected) co-oriented normal surface in
$M$. If $M= \mathbb{R}P^3 \# \mathbb{R}P^3$ then either component is
separating and therefore violates the recurrence of $M$. Hence 
$M = S^2 \times S^1$, as required. 
\end{pf}

\section{A Generalization of Novikov's Theorem}

If $M$ admits a local orientation in which each directed loop is essential,
then $\tilde M$ has no oriented loops. The results of our first section
show that any compact subset $K \subset \tilde M$ admits a transverse
measured foliation in normal form.

We now show that we can foliate $\tilde M$ globally.

\begin{thm} If $M$ admits a (not necessarily recurrent) local orientation
in which each directed loop is essential, then $\tilde M$ admits a transverse
measured foliation in normal form.
\end{thm}
\begin{pf}
The idea is to collapse the partial ordering on the vertices of the
$1$-skeleton of $\tilde M$ to a total ordering, with some kind of
geometric control, in order to construct a map $\tilde M \to \mathbb{R}$
which is an orientation-preserving embedding on each edge. 

\vskip 12pt

The proof is
by induction. At stage $i$ we will have an infinite $CW$ $2$-complex $K_i$
where $K_i$ is obtained from $K_{i-1}$ by collapsing an interval, and
where $K_0 = \tilde M^2$, the $2$-skeleton of the universal cover.
Each $2$-cell $D$ of $K_i$ will 
have the property that
the attaching map $\partial D \to K_i$ is an embedding away from possibly
finitely many points, and the induced orientation on $\partial D$ will have
exactly one maximum and one minimum. Call the two oriented subarcs of $D$
the {\it sides}. We can arbitrarily call one the left side, and one the
right side.

Let $v_l$ and $v_r$ be the two highest (with respect to the partial ordering)
vertices in $\partial D$ other than the unique maximum vertex $v$. They are
on opposite sides of $\partial D$, and are thus
incomparable in $\partial D$, but not necessarily in $K_i$. If $v_l$ and
$v_r$ are the same vertex in $K_i$, we join them by an arc in $D$ and
then collapse this arc to produce $K_{i+1}$. Otherwise, assume they are
different in $K_i$. 

It is possible that a directed path $\alpha_l$ exists from $v_r$ to $v_l$ in
$K_i^1$, or a path $\alpha_r$ from $v_l$ to $v_r$, but not both,
since $K_i^1$ contains no directed loops, and the vertices are distinct.
Suppose without loss of generality that $\alpha_l$ does not exist. 

Then if we choose a point
$p$ in the midpoint of the directed arc from $v_l$ to $v$, join $p$ to $v_r$
by an arc in $D$, and collapse this arc to produce $K_{i+1}$, the resulting
oriented $1$-skeleton will still be partially ordered.
For, if such a directed loop $\alpha$ exists, then there is a directed path
in $K_i^1$ from $v_r$ to $p$ or from $p$ to $v_r$. A directed outgoing
arc from $p$ must pass through $v$, so in the second case we would have
a directed path from $v$ to $v_r$. But there is a directed path from $v_r$
to $v$, which gives a directed loop in $K_i^1$, a contradiction.
Similarly, an incoming loop to $p$ must pass through $v_l$, so in the
first case we would have a directed path from $v_r$ to $v_l$. But by our
assumption, no such $\alpha_l$ exists. Hence $K_{i+1}^1$ is partially
ordered, and the induction step is complete. Of course, if no such $\alpha_r$
existed, we could have chosen $p$ as the midpoint of the arc from
$v_r$ to $v$.

\vskip 12pt

We see that there is a great deal of geometric control in this process:
we get to choose the $2$-cell we want to subdivide, and subdivide it
by a definite amount, dividing it into a small bigon and another cell with
one fewer vertices (of course, bigons can be completely collapsed!). 
Therefore it only takes a finite number of steps to collapse any compact
region to a stage where the edges are all of length $<$ any $\epsilon$.

We take an exhaustion of $\tilde M$ by compact sets $C_1 \subset C_2 \subset
\dots$. Then we shorten every edge of $C_1$ to length $2^{-1}$, then
every edge of $C_2$ to length $2^{-2}$, and so on until we eventually
shorten any given edge to less than $\epsilon$ in finite time. By uniformity
of this collapse we can pass to a limit. In this limit every $2$-cell has
been collapsed, and the limiting object is an infinite $1$-manifold. 
It is clear from the construction that each
edge is embedded by an orientation-preserving embedding. The associated
embeddings of the upper and lower edge of every triangle in the long edge
of the triangle are compatible across every tetrahedron, and the
induced foliation of $\tilde M^2$ therefore extends to a transverse
foliation of $\tilde M$. The local condition at vertices implies that this
foliation is non-singular. 

By construction, this foliation has no holonomy, so it admits a transverse
measure.
\end{pf}

\begin{rmk}
Suppose that $M$ is hyperbolic, and suppose that it is triangulated by
geodesic simplices. There is an $\epsilon$ and a $c$ such that if one
can show that every edge has length $\ge c$ and the angle defect at an
oriented angle (i.e. between an incoming and an outgoing edge to a vertex)
is $\le \epsilon$, then every oriented loop is homotopically essential.
To see this, observe that a piecewise
geodesic in $\mathbb{H}^2$ whose edges are all length $c$ and whose angle 
defects are all $\epsilon$ is an embedded quasigeodesic when 
$$ \frac {\pi - \epsilon} 2 > \sin^{-1}\biggl( \frac 1 {\cosh c/2}\biggr)$$
Any choice of $\epsilon, c$ satisfying the inequality above will work, by
a comparison argument.

More generally, given geometric control on $M$, it is possible in certain
circumstances to verify that every oriented loop is homotopically essential
by showing that every oriented subarc of length $<K$ for some sufficiently
large constant is quasi-geodesic with a sufficiently small coefficient of
quasi-geodicity. Since such conditions are merely sufficient but not
necessary, we do not pursue this point.
\end{rmk}

The previous theorem did not require the local orientation on $M$ to
be recurrent. However, that condition is critical for the next theorem.

\begin{thm}
Let $M$ admit a recurrent local orientation in which every oriented loop
is homotopically essential. Suppose $\tilde M$ is transversely measured
by $f:\tilde M \to \mathbb{R}$ normal on every tetrahedron, and non-singular
everywhere. Suppose $i:D \to \tilde M$ has boundary contained in a leaf.
Then $i$ can be homotoped rel. $\partial D$ to map $D$ entirely into that
leaf. 
\end{thm}
\begin{rmk}
This somewhat technical theorem is essential to what follows, and constitutes
the analogue in our context to the ``main step'' of the proof of
Novikov's theorem; see {\bf [No]}.
\end{rmk}
\begin{rmk}
Since a normal foliation is determined up to normal isotopy by its 
intersection with the $1$-skeleton, we can straighten this foliation on
each simplex in such a way that each normal disk is a flat triangle or
a quadrilateral made up of $4$ flat triangles, with respect to some
equivariant affine structure on each tetrahedron pulled back from $M$.
We assume below that this has been done.
\end{rmk}
\begin{rmk}
A posteriori, using Rosenberg's theorem, it will be shown that the disk
$D$ can be {\em isotoped} into a leaf, and not merely homotoped.
\end{rmk}
\begin{pf}
Since the homotopy property described above is open, it suffices to show that
a limit of disks which can be homotoped can itself be homotoped.

First observe that since each leaf is a normal surface, any sufficiently
small disk with boundary on a leaf can be isotoped into that leaf.

Assume that $f|_D$ has the following form.

\begin{itemize}
\item{$f \circ i(\partial D) = 1$ and $f \circ i(D) <1$ on the interior}
\item{$f$ has exactly one critical point, a minimum, on $i(D)$} 
\item{the preimages $f^{-1}(p) \cap i(D)$ foliate $i(D)$ by
concentric circles, nesting about this minimum}
\item{the minimum of $f$ on $i(D)$ is $0$}
\end{itemize}

If we can show that the theorem holds for $D$ with $f|_D$ of this form, then
we can show that it holds inductively for all $D$, by successively pushing
in innermost disks foliated as above, and reducing the number of critical
points of $f|_D$. Such a sequence of moves might involve self-intersections of
$i(D)$ with itself (i.e. it might be a homotopy rather than an isotopy) but it
will be a homotopy through immersions, since it restricts either to an isotopy
or to the identity on each piece, at each time.

We denote $$C_t = i(D) \cap f^{-1}(t), D_t = i(D) \cap f^{-1}([0,t])$$
and we let $E_t$ be the disk, a subset of $f^{-1}(t)$, whose boundary is
$C_t$, for $t<1$. The existence of $E_t$ is guaranteed by hypothesis.

If the $E_t$ lie in a compact subset of $\tilde M$, then their limit $E$
exists and is a disk, since each $E_t$ is a normal surface. Since by 
hypothesis, no such disk $E$ exists, we must suppose that the $E_t$ go off to
infinity and the leaf $E$, the subset of $f^{-1}(1)$ bounded by 
$i(\partial D)$, is noncompact. (In case $i(\partial D)$ is non-separating, 
we can consider instead
$E_t$ for $t$ extremely large, so that $E_t$ can be assumed to leave any
compact neighborhood of $i(D)$ that we choose.)

We establish the following lemma

\begin{lem}
In the above context there is a $K$ such that for every normal subdisk $D$
of a leaf whose boundary bounds another normal disk $D'$ in $M$, such that
$D$ and $D'$ are homologous, if $|D|$ denotes the number of normal disks
in $D$, we have an inequality
$$|D| \le K|D'|$$
\end{lem}
\begin{pf}
Let $\alpha$ be a loop in the $1$-skeleton of $M$ which passes through every
edge. We can find such an $\alpha$ by the assumption that the $1$-skeleton
of $M$ is recurrent.
Since $D,D'$ were assumed to be homologous, the surface $D \cup D'$ is
null-homologous in $\tilde M$ and it therefore projects to a null-homologous
surface in $M$, whose algebraic intersection number with $\alpha$ is therefore
zero. Now, there is some constant $c_1>0$ such that the geometric intersection
$|\alpha \cap D'|\le c_1|D'|$, since $\alpha$ passes only finitely many
times through each edge in $M^1$. Moreover, since $D$ is a subdisk of a leaf,
the geometric intersection of $\alpha$ with $D$ is equal to their algebraic
intersection, which is at least $c_2|D|$ for some $c_2>0$ depending on the
maximal order of an edge in $M^1$. From $c_1$ and $c_2$ we can find our
constant $K$.
\end{pf}

Now, each leaf $E_t$ in homologous to some subdisk of $i(D)$. By simplicial
approximation, these subdisks of $i(D)$ are all approximated by normal disks
of bounded size. By our lemma, therefore, there is a bound on the number of
normal disks in $E_t$, contradicting our assumption that the $E_t$ went off
to infinity. This contradiction establishes the theorem.

\end{pf}

\begin{cor}
Every leaf in $\tilde M$ as above is incompressible.
\end{cor}
\begin{pf}
This is immediate by the loop theorem.
\end{pf}

\begin{cor}
$M$ as above is irreducible or $S^2 \times S^1$.
\end{cor}
\begin{pf}
If $M \ne S^2 \times S^1$ and it is reducible, then there is an embedded
separating sphere in $M$. This lifts to an embedded sphere in $\tilde M$.
Then our result follows word for word the 
proof of Rosenberg's Theorem {\bf [Ro]},
after remarking that every leaf in $\tilde M$ is a possibly disjoint union
of $\mathbb{R}^2$'s and $S^2$'s. If any leaf contains an $S^2$, the Reeb 
stability theorem, together with our result, shows that $\tilde M$ is
$S^2 \times \mathbb{R}$ (since, inductively, no ``first'' leaf can become
non-compact).
Since $M$ is compact and orientable, it is
$S^2 \times S^1$ or $\mathbb{R}P^3 \# \mathbb{R}P^3$. But if 
$M = \mathbb{R}P^3 \# \mathbb{R}P^3$ then as in the previous section
the image of an $S^2$ separates $M$ 
and contradicts recurrence. If $\tilde M$ is foliated by $\mathbb{R}^2$'s
it is irreducible.
\end{pf}

\begin{cor}
If $\tilde M$ is $\delta$-hyperbolic with respect to (any) complete metric 
pulled back from $M$, then each leaf in $\tilde M$ foliated as above by
incompressible leaves is $\delta'$-hyperbolic, where $\delta'$ depends only
on $\delta$ and the combinatorics of $M$.
\end{cor}
\begin{pf}
From the proof and the statement of theorem 5.2 we can deduce that for every
normal subdisk $D$ of a leaf whose boundary bounds another normal disk $D'$
in $\tilde M$, we have an inequality $|D| \le K |D'|$ for some $K$ 
depending only on the triangulation of $M$. For, 
$\pi_2(\tilde M) = H_2(\tilde M) = 1$, and therefore any two disks in
$\tilde M$ with the same boundary are homologous.

Now, we know that $\delta$-hyperbolic Riemannian manifolds are characterized by
the fact that minimal spanning disks satisfy linear isoperimetric
inequalities (see, for instance, {\bf [Gr]}). Our comparison lemma allows
us to establish a similar isoperimetric inequality for subdisks of a leaf,
and therefore each leaf is also $\delta'$-hyperbolic. One can see that
$\delta'$ depends only on the combinatorics of $M$ and on $\delta$.
\end{pf}

\begin{rmk}
Note that Rosenberg's theorem shows that disks can be isotoped into leaves,
rel. boundary. For if $D_1 \simeq D_2$ rel. boundary, we can consider
intersections $D_1 \cap D_2$, which will be a collection of circles. But
innermost circles in the intersection will define $S^2$'s in $D_1 \cup D_2$
which must bound $B^3$'s. Pushing $D_1$ across this $B^3$ we can isotope
it to reduce the number of components of $D_1 \cap D_2$. Inductively, this
shows we can isotope $D_1$ to $D_2$ rel. boundary.
\end{rmk}

\begin{rmk}
Notice that the hypothesis that every oriented loop in $M$ be homotopically
essential is used {\it only} to establish the existence of a foliation on
$\tilde M$ in normal form, transverse to the orientations pulled back from
$M$. If we are given this foliation on $\tilde M$ as a {\it hypothesis}, say
if the orientation on $M$ was inherited from a foliation in the first place,
then the proof of the theorem still goes through, and we can show that 
every oriented loop in $M$ is homotopically essential as a {\it consequence}
of the existence of the foliation on $\tilde M$. Thus, our technique gives
a new proof (conceptually similar, though perhaps technically easier) of
the theorem of Novikov that circles transverse to taut foliations are 
homotopically essential.

More explicitly, given a taut foliation on $M$, we can take a sufficiently
fine triangulation and orientations on the edges to obtain a local orientation.
Lifting this foliation to the universal cover, our theorem applies to show
that every leaf in $\tilde M$ is incompressible and is therefore a disk
or a sphere. Such a foliation has no holonomy, so it admits a transverse
measure. Since we are on $\tilde M$, a transverse measure is given by
integrating an exact $1$-form $df$. A homotopically inessential 
transverse circle
in $M$ lifts to a transverse circle in $\tilde M$. But $df$ is positive 
everywhere on the tangent vector to this circle which is absurd.
\end{rmk}

\begin{rmk}
Our result shows that the foliations on the universal cover are, roughly
speaking, minimal surfaces with respect to the weights on edges determined
by an oriented cycle in $M$ passing through every edge. This is in some sense
a combinatorial (non-deterministic) volume-preserving flow on $M$. This flow
is used by Gabai in {\bf [Ga2]} to prove a refinement of a theorem of
Roussarie and Thurston: if $S$ is an immersed incompressible surface in
$M$ admitting a taut foliation, then $S$ can either be homotoped into a
leaf, or can be homotoped to have only saddle-type tangencies with the
foliation. 
\end{rmk}

We can produce a qualitative topological refinement of theorem 5.2 by a 
more delicate argument.

\begin{thm}
Suppose $M$ is not prime, and suppose its $1$-skeleton is ordered in such
a way that all the conditions of the above theorem except recurrence are
satisfied. Then there exists an embedded null-homologous co-oriented
(not necessarily connected) normal surface in $M$ such that every transverse 
arc is outgoing.
\end{thm}
\begin{pf}
As in theorem 5.1 we can produce a measured foliation on $\tilde M$.
If every leaf were incompressible, $M$ would be prime. Therefore some leaf
is compressible, and by the loop theorem, there is an embedded disk $D$ in
$\tilde M$ whose boundary lies on a leaf, such that there is no disk in that
leaf with the same boundary. We adopt the notation from the proof of 
theorem 5.2. and assume the restriction of the foliation to
$D$ is by concentric circles, by induction.

\vskip 12pt

Recall that $D$ is a disk whose boundary lies in a leaf of the foliation
of $\tilde M$, that $E$ is the subset of this leaf bounded by $\partial D$,
and that $E$ is a ``limit'' of the $E_t$ - disks contained in the foliation
whose boundaries are circles making up a concentric foliation of $D$.

If $R$ is the region in $\tilde M$ bounded by $D$ and $E$ then every
increasing arc that passes through $E$ must leave $R$. Now, $E$ is made up
of normal triangles and quadrilaterals, so there are outgoing arcs within
bounded distance of any point on $E$. Since the only incoming arcs to $R$
pass through $D$, and since every vertex is the endpoint of some arc, 
there are arbitrarily long
paths contained in $\tilde M^1$ passing through $D$ and 
contained entirely within $R$. Since each sufficiently long arc eventually 
passes through
each $E_t$, and since each $E_t$ bounds a compact region of $R$, we can
extract a subsequence of these paths which converge on compact sets (since
they are all simplicial!), and find an infinite increasing path 
$\gamma \in \tilde M^1$ whose initial point is in $D$ and which is
contained entirely inside $R$. Notice that $\gamma$ intersects each $E_t$
exactly once for all $t>t_0$, the value of $f$ at the initial point of
$\gamma$. Let $\gamma_t$ denote the initial segment of $\gamma$ from
$0$ to $t$.

Since the $D_t$ converge to $D$, there is some $t$ after which all the
$C_t$ are normally isotopic to $\partial D$. If we truncate $R$ by only 
considering the region above this $t$, we can replace $D$ by a slight
perturbation of $E_t$, and therefore we can assume, without loss of
generality, that every increasing simplicial arc passing through $D$ is 
incoming to $R$. 

By the compactness of $M$, we can assume there is some $p \in M^0$ such
that $\gamma$ passes through infinitely many lifts of $p$ which we call
$p_0,p_1, \dots$. Let $\alpha_i \in \pi_1(M)$ be such that
$\alpha_i(p_i) = p_0$. After passing to a subsequence and re-ordering if
necessary, we can assume that the translates $\alpha_i(D)$ are disjoint
and the collection is embedded.

\vskip 12pt

Since our foliation of $\tilde M$ is not necessarily $\pi_1(M)$ equivariant,
the image of subsets of the leaves can intersect. We need to investigate
these intersections more closely.
We have the following lemma, which controls the orientations on 
$\alpha_i(E_t),\alpha_j(E_s)$ when they intersect in an isolated point of
tangency.

\begin{lem}
If some normal subsurfaces of $\alpha_i(E_t),\alpha_j(E_s)$ intersect
in an isolated point of tangency for some $i,j,s,t$ then the transverse
orientations to $\alpha_i(E_t),\alpha_j(E_s)$ agree at this point.
\end{lem}
\begin{pf}
If the intersection is at a vertex, then a neighborhood of the intersection
in either surface separates the star of the vertex. By the definition of 
the foliation on $\tilde M$, the outgoing and the incoming edges to the
vertex lie in different components of the star. Since both these collections
are non-empty, there is a monotone arc transverse to both surfaces at
the point of intersection, whose orientation agrees with the transverse
orientations on both surfaces at this point.

If the intersection is in some tetrahedron $\Delta$, then there are a pair
of normal disks in $\alpha_i(E_t)\cap \Delta, \alpha_j(E_s)\cap \Delta$
which intersect in an isolated point of tangency. Since any two normal
disks in a tetrahedron intersect a common edge, this edge is transverse
to both surfaces at this point and their transverse orientations therefore
agree.
\end{pf}

Since each $\alpha_i(\gamma_t)$ must pass through $p_0$ for some $t$, and
since $\alpha_i(\gamma_t)$ lies outside $R$ for sufficiently small $t$,
there is some increasing sequence $t_i$ such that $\alpha_i(\gamma_{t_i})$ 
intersects $D$. 
As we increase $t$ past $t_i$, the surface $\alpha_i(E_t)$ intersects a 
collar neighborhood of $D$ in an annulus insulating 
$\alpha_i(\gamma)$ from $\partial D$, and expanding concentrically with
$t$. Let $B_{i,0}(t)$ be the circle of intersection with $D_t$ for $t$ near
$t_i$. Then as $t$ increases, we may push part of $B_{i,0}(t)$ over the edge
of $D_t$ and up into $E_t$. That is, we think of $B_{i,0}(t)$ as the 
appropriate component of $\alpha_i(E_t)\cap (E_t \cup D_t)$.

\begin{lem}
There is some finite $n$ such that there are at most $n$ embedded circles
and properly embedded arcs made from pieces of $\alpha_i(E_t)\cap D$ 
extending normally to a simplicial collar of $D$ such that no two of the 
collection are normally isotopic in a neighborhood of $D$.
\end{lem}
\begin{pf}
These circles and arc are the meridians of embedded normal annuli and
arcs $\times I$ transverse to $D$. Therefore they bound embedded subdisks
of $D$, and since they are normal, the length of the circles and arcs is
bounded. Hence each annulus and arc $\times I$ is composed of a bounded
number of pieces, and since the simplicial neighborhood is finite,
there are only finitely many normal surface types represented by them.
\end{pf}

Notice that since $t_j<t_i$ for $j<i$, the surface made up from
$\alpha_j(E_t)\cup \alpha_j(D)$ for appropriate $t$ 
separates $\alpha_i(D)$ from $p_0$ for
$i>j$. In particular, $\alpha_i(\gamma_t)$ must pass through each
$\alpha_j(D)$ for $j<i$ before passing through $D$. Suppose 
$\alpha_i(\gamma_t)$ exits $\alpha_j(E_t)$ for some $j<i$ and $t<t_i$. Then
again it can never reach $p_0$. We denote the circles of intersection of
the $\alpha_i(E_t)$ with $\alpha_j(E_t \cup D_t)$ by $B_{i,j}(t)$. Notice
that $$B_{i,i}(t) = \alpha_i(\partial D_t)$$
which is normally isotopic to $\alpha_i(\partial D)$ for large $t$.

The $B_{i,j}(t)$, with $i\le k$ divide $\alpha_j(D)$ for fixed $j$ into a 
collection of regions. By our previous comment, each $B_{i,j}(t)$ bounds
the region containing $\alpha_k(\gamma) \cap \alpha_j(D)$ for $i\le k$.
Let $B^k_j(t)$ be the boundary of the subregion of $D$ containing
$\alpha_k(\gamma)\cap \alpha_j(D)$. This is a circle contained entirely 
within $\alpha_j(D)$. Call $B^k_j(t)$ an {\em innermost} circle.

\vskip 12pt

Since there are only finitely many possibilities for the $B^k_j(t)$ for each
$j$, up to normal isotopy, by choosing $i$ very large we can find 
$B^i_j(t)$ such that 
$\alpha_k\circ \alpha_j^{-1}(B^i_j(t))$ is normally isotopic to $B^i_k(t)$.
We know that the innermost circle $B'$
made up of $B_{l,k}(t)$ with $j \le l \le i$ lies outside $B^i_k(t)$, since it
is the innermost of fewer circles. By definition, there is an annulus 
made up of pieces of $\alpha_i(E_t)$ interpolating between $B^i_j(t)$ and
$B'$. This annulus is embedded, since it bounds some image of $\gamma$, and
we can take an innermost such. In more detail, this annulus is the boundary
of the connected region in the complement of the relevant $\alpha_i(E_t)$
which bounds the relevant image of $\gamma$. Since each $\alpha_i(E_t)$ bounds
this image of $\gamma$, such a region exists. Our orientation lemma implies
it is an annulus.

\vskip 12pt

Since each $B_{i,k}(t)$ bounds some subdisk of $\alpha_i(E_t)$, we can
cut and paste an innermost disk which bounds $B'$. Together with the
subdisk $D'$ of $\alpha_j(D_t)$ bounded by $B^i_j(t)$ this gives a 
(topological) sphere $S$ bounding a $B^3$ in $\tilde M$ such that $B^i_k(t)$ 
and the subdisk $\alpha_k \circ \alpha_j^{-1}(D')$ 
of $\alpha_k(D_t)$ that it bounds is entirely contained inside the region
bounded by $S$. One should be careful to note that the sphere in question
bounds a $B^3$ because it is contained in the region $R$ whose interior is 
foliated with disks, and therefore irreducible.

\vskip 12pt
Let $\beta= \alpha_k \circ \alpha_j^{-1}$, and let $N = \tilde M/<\beta>$.
Then since $S$ bounds a ball in $\tilde M$, its image under the projection
to $N$ is some compact submanifold on $N$. Its boundary cannot contain
any piece of $D'$, since $D'$ is interior to some translate of the ball
bounded by $S$. If $N$ is non-compact, this boundary is non-empty and
by construction is a null-homologous normal surface in $N$ made up entirely of
projections of pieces of $E_t$. Call this surface $H$ and consider its
projection $\pi(H)$ to $M$.

If the boundary is empty, then $N$ is compact, and $M$ is $S^2 \times S^1$
by an earlier result. 

\vskip 12pt

Since $\pi(H) \subset M$ is made up entirely of pieces of the image of
$E_t$ under the projection $\tilde M \to M$, it is represented by an
embedded (possibly disconnected) normal surface $G = G_1 \cup G_2 \cup \dots
\cup G_n$. By our orientation lemma, each $G_i$ is co-oriented such that
every transverse arc in the $1$-skeleton is outgoing.
\end{pf}

\begin{rmk}
This theorem slightly weakens the condition of recurrence to prove 
irreducibility.

Notice that such a surface $G$ is a finite Haken sum of fundamental
normal surfaces co-oriented compatibly with the orientation on the 
$1$-skeleton. One can check algorithmically whether some $\mathbb{Z}$-linear
combination of such fundamental surfaces can be trivial in $H_2$.
\end{rmk}

\section{Is Homotopically Essential Essential?}

In this section we show that, at least to prove irreducibility, the condition
that every oriented loop in the $1$-skeleton be homotopically essential can
be substantially weakened, and even weakened to an easily checkable condition.

\begin{defn}
Let $M$ be a triangulated $3$-manifold. For any $m \in \mathbb{Z}$,
the simple combinatorial $m$-germ at a vertex $p$, denoted $\tilde M_m(p)$,
 is the simplicial complex
obtained in the following way:

\begin{itemize}
\item{Let $N_m(p)$ be the disjoint union of the simplicial neighborhoods of
simplicial paths of length $m$ in $M$ with initial vertex $p$.}
\item{Obtain $\tilde M_m(p)$ as the quotient space of $N_m(p)$ by identifying
endpoints of two distinct paths which have the same endpoint in $M$, and
which bound an (immersed) simplicial disk in $M$ of simplicial area $\le m$.}
\end{itemize}

\end{defn}

\begin{rmk}
It is clear from the definition that the complex $\tilde M_m(p)$ can be
algorithmically constructed. Note that we could fine-tune the relative
sizes of paths and disks in the definition to more accurately capture 
approximations to the germ of $\tilde M$ at $p$ using estimates of
an isoperimetric inequality for $M$. 
\end{rmk}

\begin{thm}
Let $M$ admit a local orientation such that for an appropriate, explicitly
computable constant $k$ depending only on the triangulation and the orientation
of $M$, the complex $\tilde M_k(p)$ with the induced orientations on the
$1$-skeleton has no oriented loops. Then $\pi_2(M)=1$ or $M=S^2\times S^1$.
\end{thm}
\begin{pf}
If $M$ is not $S^2\times S^1$ and is reducible, then there exists a
separating normal $S^2$. This $S^2$ lifts to $\tilde M_k(p)$ for sufficiently
large $k$. If $\tilde M_k(p)$ has no oriented loops, it admits a transverse
co-oriented foliation in normal form.

As in the proof of our earlier theorem, we consider the intersection of
the $S^2$ with this foliation, and start to push innermost disks into the
leaves. Our earlier estimates for the simplicial size of these disks
still holds, since $M^1$ is recurrent. Every subdisk of the homotopy is
within an (easily computable) distance from some fixed $p \in S^2$, and all
the disks that we are pushing are of simplicial size bounded by some
computable constant times the simplicial size of the $S^2$, so this
homotopy can be carried out within $\tilde M_k(p)$ (i.e. we never push
over the boundary).

If this $S^2$ can be pushed entirely into an $S^2$ leaf, then there is a
separating $S^2$ in $M$ oriented compatibly with the $1$-skeleton, which
contradicts recurrence. Otherwise, the $S^2$ can be pushed entirely into a
disk, and therefore was null-homotopic in $M$.

The theorem is proved once we observe that we can bound the simplicial size
of the smallest separating embedded normal homotopically essential $S^2$ in 
terms of the triangulation of $M$. 
\end{pf}

\section{Questions}

It is natural to ask to what extent some of the technical hypotheses in this
paper can be removed. In particular, the following questions seem outstanding:

\begin{enumerate}
\item{Is there an algorithm to check whether every oriented loop in a
local orientation is homotopically essential?} 
\item{To what extent is our combinatorial structure weaker than the existence
of a taut foliation in normal form?}

\begin{rmk}
The referee has pointed out that good candidates for local orientations 
with homotopically essential loops which
do not admit transverse foliations might be found by investigating certain
graph manifolds, in particular those obtained from products 
$(\text{punctured surface}) \times S^1$ by appropriate glueings along the
boundary tori. Work of Brittenham, Naimi and Roberts {\bf [BNR]} shows that
many such graph manifolds do not admit any taut foliations whatsoever. On
the other hand, such manifolds certainly admit local orientations, and one
expects that the condition that oriented loops be homotopically essential
can be satisfied in many cases.
\end{rmk}
\item{If a triangulation of $M$ admits a local orientation in which every 
oriented loop is homotopically essential, to what extent can this 
combinatorial structure be extended over a refinement of the triangulation?}
\item{Is there some (computable) bound on the number of subdivisions of a
triangulation necessary to put it in normal form with respect to an existing
($C^1$) foliation on $M$?}
\item{Can the finiteness of the combinatorial structure be used to advantage
in addressing questions of the virtual existence of such a structure - i.e.,
when does there exist a finite cover of $M$ which admits a local orientation
in which each loop is homotopically essential? or lies in an open half-space
of $H_1$?}
\begin{rmk}
This question is intimately related to Thurston's famous conjecture that
every hyperbolic $3$-manifold has a finite cover which fibers over the
circle. In fact, this question was our main original motivation for studying
the interaction of foliations with finite combinatorial structures. 
\end{rmk}

\item{If $\alpha$ is an oriented embedded loop in $M^1$ which has a local
orientation in which each oriented loop is homotopically essential, is it
true that all Dehn surgeries on $\alpha$ with slope sufficiently close to
$(1,0)$ give manifolds whose $1$-skeleton can be similarly oriented without 
changing the triangulation or the orientation on $M - {\rm nbhd}(\alpha)$?}
\item{To what extent can the geometry of leaves in $\tilde M$ be controlled?}
\item{Can one extend the results of the last few sections to triangulations
with local orientations on subgraphs of the $1$-skeleton?}
\begin{rmk}
Our local orientations are dual to branched surfaces such that complementary
domains are sutured balls. More generally, branched surfaces whose 
complementary domains are sutured manifolds admitting taut foliations 
can (with certain technical hypotheses) carry only incompressible surfaces
(see {\bf [Oe]}). It seems plausible that if one controls the complementary
regions, homotopically essential local orientations on subsets of 
$M^1$ might have many nice properties.
\end{rmk}
\end{enumerate}

\section{References}

\begin{enumerate}
\item{{\bf [Be]} M-T. Benameur, Triangulations and the Stability Theorem for
Foliations, {\it Pacific J. Math.}, {\bf 179} (1997) pp.221-239.}
\item{{\bf [BNR]} M. Brittenham, R. Naimi and R. Roberts, Graph Manifolds
and Taut Foliations, {\it J. Diff. Geom.}, {\bf 45} (1997) pp.446-470.}
\item{{\bf [Ga]} D. Gabai, Problems in Foliations and Laminations.
in ``Geometric Topology'', (edited by W. Kazez), proceedings of the 1993
Georgia International Topology Conference, Vol. 2, part 2, pp. 1-33.}
\item{{\bf [Ga2]} D. Gabai, A Combinatorial Volume Preserving Flow on
Taut Foliations. {\it preprint}}
\item{{\bf [Gr]} M. Gromov, Hyperbolic Groups. in ``Essays in Group
Theory''. MSRI Publ. {\bf 8} (1987) pp. 75-263.}
\item{{\bf [No]} S. Novikov, Topology of Foliations, {\it Trans. Moscow
Math. Soc.}, {\bf 14} (1965) pp. 268-305.}
\item{{\bf [Oe]} U. Oertel, Homology Branched Surfaces: Thurston's Norm
on $H_2(M^3)$. in ``Low-Dimensional Topology and Kleinian Groups''.
LMS Lecture Note Series 112. (1986) pp. 253-272.}
\item{{\bf [Ro]} H. Rosenberg, Foliations by Planes, {\it Topology} 
{\bf 7}, (1968), pp. 131-138.}
\item{{\bf [Su]} D. Sullivan, Cycles for the Dynamical Study of
Foliated Manifolds and Complex Manifolds, {\it Inventiones Math.}, 
{\bf 36} (1976) pp. 225-255.}
\end{enumerate}

\end{document}